\documentclass[12pt]{article}
\usepackage[explicit]{titlesec}
\usepackage[none]{hyphenat}
\usepackage{amsmath,amssymb,amsfonts,mathrsfs,bm}

\usepackage{soul}
\setul{.6pt}{.4pt}

\usepackage{geometry}
\usepackage{fancyhdr}
\usepackage[labelfont={footnotesize,bf} , textfont=footnotesize]{caption}
\titleformat{\section}
{\normalfont\bfseries}{\thesection.}{0.5em}{\MakeUppercase{\textbf{#1}}}
\titlespacing\section{0pt}{0pt}{-10pt}
\titleformat{\subsection}
{\normalfont\itshape}{\thesubsection.}{0.5em}{\textit{#1}}
\titlespacing\subsection{0pt}{0pt}{-8pt}

\makeatletter
\newcommand\sixteen{\@setfontsize\sixteen{17pt}{6}}
\renewcommand{\maketitle}{\bgroup\setlength{\parindent}{0pt}
	\begin{flushleft}
		\sixteen\bfseries \@title
		\medskip
	\end{flushleft}
	\textit{\@author}
	\egroup}
\makeatother

\usepackage{hyperref}

\usepackage{caption}
\usepackage{graphicx,subfigure}
\graphicspath{{figures/}}

\usepackage{booktabs}
\usepackage{multirow}
\usepackage{array}

\usepackage[american]{babel}
\usepackage{microtype} 
\usepackage{enumerate}

\usepackage{tikz}

\title{Fatigue crack propagation in carbon steel using RVE based model}

\author{
	Zhenxing Cheng$^{a}$, Hu Wang$^{a}$\footnote{Corresponding author. Email: wanghu@hnu.edu.cn}, Gui-Rong Liu$^{b}$ \medskip \\
	$^{a}$State Key Laboratory of Advanced Design and Manufacturing for Vehicle Body, Hunan University, Changsha, 410082, PR China \\
	$^{b}$Department of Aerospace Engineering and Engineering Mechanics, University of Cincinnati, Cincinnati, Ohio, 45221, United States\\}

\begin{document}
	
	\vspace*{.01 in}
	\maketitle
	\vspace{.12 in}
	
	\section*{abstract}
	A representative volume element (RVE) based multi-scale method is proposed to investigate the mechanism of fatigue crack propagation by the molecular dynamics (MD) and the  extended finite element methods(XFEM) in this study. An atomic model of carbon steel plate is built to study the behavior of fatigue crack at the micro scale by MD method. Then the RVE model  for fatigue crack propagation should be built by  fitting the data which was obtained from the MD result with the Paris law model. Moreover, the effect of micro-structural defects including interstitial atoms, vacancies have also been considered in this study. The results indicate that the micro-structural defects can deeply influence the values of Paris law constants and the life of the specimen can be evaluated by the proposed method.
	
	\textit{Keywords}: Molecular dynamics, Fatigue crack, Extended finite element method,  Multi-scale
	
	\vspace{.12 in}
	
	
	\section{introduction}
	
	Detection and simulation of failures are vital to ensure the stability of structures in engineering fields, e.g., aircraft, watercraft, vehicle. Fatigue fracture, as a main part of structural failures, usually  happens at load levels below the yield stress and eventually developed to unstable fracture. Therefore, it is  important to simulate the propagation of fatigue cracks and make predictions of life cycles before the destructive fracture occurs. In recent decades, extensive studies have been published to investigate the mechanisms of fatigue crack in metallic materials \cite{Schutz1996,Forman1967,Maierhofer2014,Chudnovsky2014,Correia2016,Borges2020,Duan2020}. Several kinds of functions or laws have been proposed to describe and predict the rule of fatigue crack propagation including polynomial \cite{Davies1973} , exponential \cite{Frost1958,Zheng2005}, power functions \cite{DeIorio2012}.  A popular model for fatigue crack propagation is the Paris law in which the relationship between cyclic crack growth rate $da/dN$ and the stress intensity factor $\Delta K$ is  described by  a power function \cite{schreiber2020phase}. In the early 1960s, Paris et al. proposed a effective semi-empirical model, known as Paris Law later, to fit the curve between $da/dN$ and $\Delta K$ \cite{Paris1963}. After that, plenty of studies have been proposed to modify the Paris law for a wider range of applications by considering more parameters \cite{Bazant1993,Kirane2016,Kim2000,Zheng1999}. However, those studies mentioned above are based on the macroscopic research while the micro-structure is not considered, i.e., defects in the specimen which can influence parameters of the material and crack are ignored, e.g., dislocations, vacancies, slip bands, twins, voids, inclusions. Therefore, it is vital to accurately evaluate whether the defects make an effect on fatigue crack propagation or not. In view of these facts, it is necessary to understand the fatigue crack propagation on the micro-scale. Molecular dynamics (MD) method,  as a popular numerical method for detailed microscopic modeling on the  molecular scale,has been widely used to study the behavior of fatigue crack in different types of crystalline materials recently \cite{huang2012discrete,leung2014atomistic,zhang2017mechanisms,Yasbolaghi2020}. Tang and Horstemeyer et al. employed MD method to study the fatigue crack propagation behavior of single magnesium crystal with different orientations \cite{tang2010fatigue} and then reviewed the studies about atomic simulations of fatigue crack propagation in nickel and cooper \cite{horstemeyer2010nanostructurally}. Ma et al. simulated the fatigue crack propagation in iron by MD method and calculated the Paris law constants with the MD simulation \cite{ma2014molecular}. After that, cohesive models based MD simulations were used to investigated the behavior of fatigue crack propagation in both single and poly crystals \cite{wu2015molecular,Lu2020,Li2015}. Three dimensional fatigue cracks were also investigated to reveal the manner of crack evolution by  MD based atomic simulations \cite{Uhnakova2011,Uhnakova2012}.
	
	To our knowledge, few works have investigated the behavior of fatigue crack propagation on different length scales, e.g., macro-scale, micro-scale, nano-scale (or atomic scale).  Therefore, a representative volume element (RVE) based multi-scale method is proposed to study the mechanism of crack propagation in this study, in which the extended finite element (XFEM) \cite{Moes1999,Belytschko2009} and MD methods are used to study the behavior of macro and micro crack propagation respectively.  Moreover, the effect of voids and inclusions  have also been considered in this study.
	
	The remainder of this study is as follows. The basic theory of Paris law and MD  are introduced in section 2. The description of atomic model and RVE based multi-scale method  can be  found in section 3. Results and discussions can be found in section 4. Finally, some conclusions of this study are given in section 5.


	\section{Basic theories}
	
	\subsection{Paris law  for fatigue crack propagation}
	
\begin{figure*}[htbp] 
		\centering
		\includegraphics[]{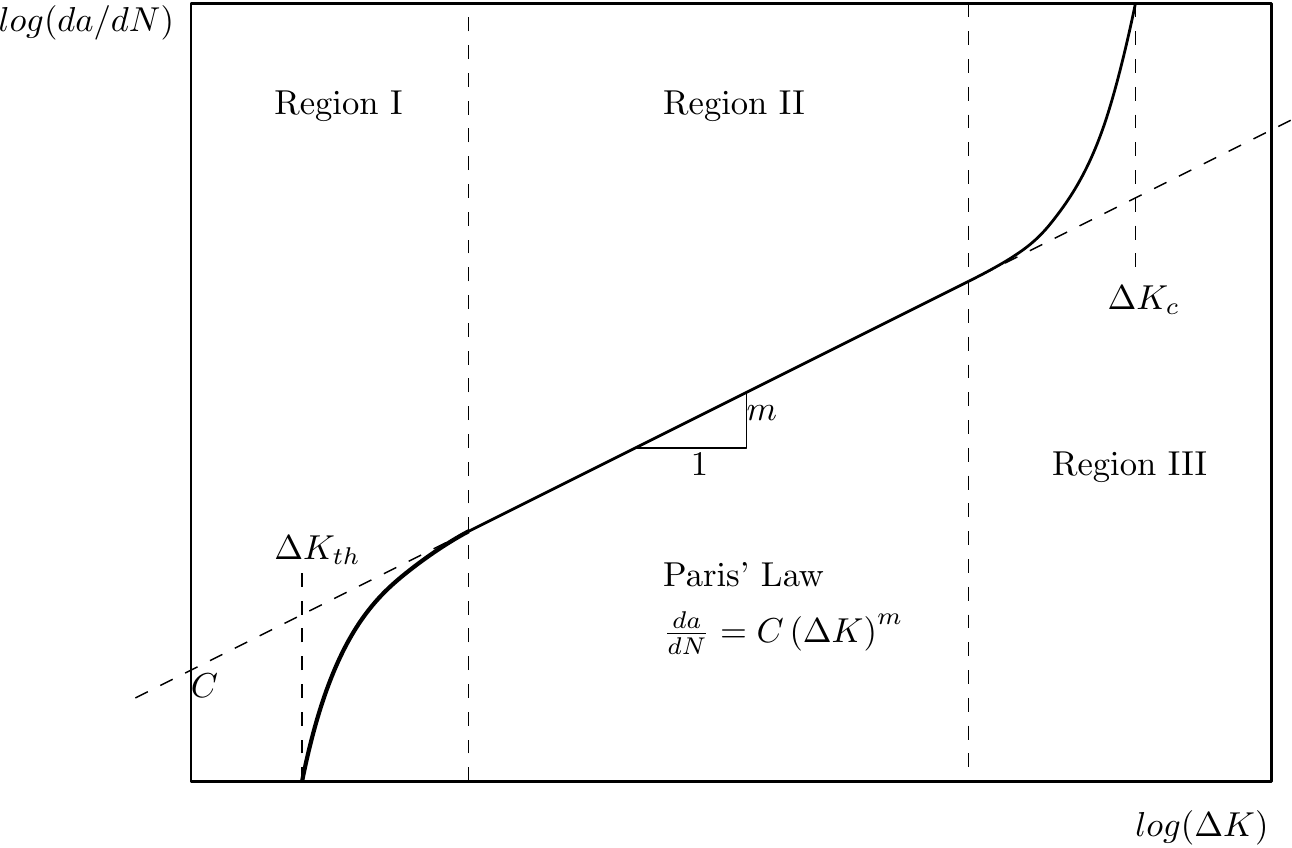}
		\caption{Fatigue growth rate curve.}
		\label{fig:paris}
	\end{figure*}

	A very common and popular model for fatigue crack propagation is the Paris model which gives the mathematical formulas for stable propagation of fatigue crack. The Paris model is described by Eq. \eqref{eq:1}, where the symbol $da/dN$ means the increment of crack during each loading cycle, and $\Delta K$ is the amplitude of stress intensity factor around the crack tips. Figure \ref{fig:paris} shows the typical fatigue growth rate curve by log-log plot which is also known as $da/dN$ versus $\Delta K$ curve. The curve can be divided into three parts: regions \MakeUppercase{\romannumeral1, \romannumeral2} and \MakeUppercase{\romannumeral3} \cite{Paris1963}. Region \MakeUppercase{\romannumeral1} represents the early stage of a fatigue crack and $\Delta K_{th}$ is the crack growth rate where the long fatigue crack starts growing. Region \MakeUppercase{\romannumeral2} represents  the stable  crack growth zone where the Paris law relationship can be described by a linear equation (Eq. \eqref{eq:2}) \cite{Bozic2011}. The mathematical equation contains two parameters $m$ and $C$, where $m$ means the slope of the line in Fig. \ref{fig:paris} and $C$ is the y-axis-intercept. The constants $m$ and $C$ are usually determined in experiments \cite{Carrascal2014,Ancona2016,Chauhan2016,Branco2009,Branco2012}, but since the experimental determination is tedious and time consuming, some numerical methods are also used to calculate these constants recently \cite{Yasbolaghi2020,ma2014molecular,Horstemeyer2010,Mlikota2017}.
	
	\begin{equation} \label{eq:1}
		\frac{da}{dN}=C\left ( \Delta K\right )^{m}.
	\end{equation}

	\begin{equation} \label{eq:2}
		log(\frac{da}{dN})=m(log\left ( \Delta K\right ))+log(C).
	\end{equation}
	
	\subsection{Molecular dynamics method}
	
	Molecular dynamics is a numerical method which is used to computationally simulate the time evolution of  atoms or molecules by Newton's equations of motion, i.e., a typical MD simulation can be described by the position and momentum of each atom or molecule. The dynamics of atoms can be described as Eq. \eqref{eq:3},
	
	\begin{equation} \label{eq:3}
	m_i\bm{a}_i = \bm{f}_i = -\triangledown U({\bm{r}_i}),
	\end{equation}
	
	\noindent where $i$ represents the number of atom, and $m_i$, $\bm{a}_i$, $\bm f_i$, $\bm r_i$ denote the mass, acceleration, inter-atomic force, position of atom $i$ respectively. The symbol $U$ denotes the inter-atomic potential which can be obtained by experiments. A typical inter-atomic potential file calculated by Embedded Atom Method (EAM)  is used in this study\cite{daw1983semiempirical}. Generally, the inter-atomic force can be calculated by Eqs. \eqref{eq:4}and \eqref{eq:5},  where $U$ is a function of the energy of each atom ($U_{ij}$), which relatives with positions of atoms $i$ and $j$. The symbol $r_{cut}$ is the cut-off radius which means the threshold distance that atoms do not interact directly.  An EAM potential file proposed by Hepburn and Ackland is used in this study \cite{hepburn2008metallic,becker2013considerations,hale2018evaluating}.
	
	\begin{gather}
	U(\bm r_i)=\sum_i U_{ij} (\bm r_i,\bm r_j), \label{eq:4}\\
	\bm f_i = -\triangledown U({\bm r_i}) = 0, \quad if \quad r_{ij}>r_{cut}. \label{eq:5}
	\end{gather}
	
	Usually, a large number of motion equations need to be solved during the MD simulation. Therefore, the Velocity-Verlet algorithm is used to solve the motion equations with considerable accuracy \cite{omelyan2002optimized}. However, only positions and velocities can be obtained by the Velocity-Verlet algorithm, so Swenson et al. suggested a definition of virial stress to calculate the atomic stress \cite{swenson1983comments}. Atomic scale virial stresses are equivalent to the continuum Cauchy stresses \cite{subramaniyan2008continuum}. The stress contains two parts, potential and kinetic energy parts, which is defined as
	
	\begin{equation} \label{eq:6}
	\sigma_{xy}=\frac{1}{V^i}\sum_i \left[ \frac{1}{2}\sum_{j = 1}^N \left( \bm r_x^j - \bm r_x^i \right) \bm f_y^{ij}-m^i \bm v_x^i \bm v_y^i \right] ,
	\end{equation}
	
	\noindent where the subscripts $x$ and $y$ represent the Cartesian components and $V^i$ means the volume of the atom $i$.  Other symbols are described above. Specially, the symbol $f_y^{ij}$ is the $y$ direction force on atom $i$ induced by atom $j$, $v_x^i$, $r_x^i$ are the velocities and relative position of atom $i$ along the $x$ direction. Then the Von Mises stress $\bar{\sigma}$ can be calculated by
	
		\begin{equation}
		\bar{\sigma}=\frac{1}{\sqrt{2}}\sqrt{(\sigma_x-\sigma_y)^2+(\sigma_y-\sigma_z)^2+(\sigma_z-\sigma_x)^2+6(\tau_{xy}^2+\tau_{yz}^2+\tau_{zx}^2)},
		\end{equation}
	
	\noindent where $\sigma_x, \sigma_y, \sigma_z$ are the normal stresses and $\tau_{xy}, \tau_{yz}, \tau_{zx}$ are the tangential stresses.

	\section{Methodology}

	\subsection{Framework of the RVE based multi-scale method}
\begin{figure*}[htbp] 
	\centering
	\includegraphics[]{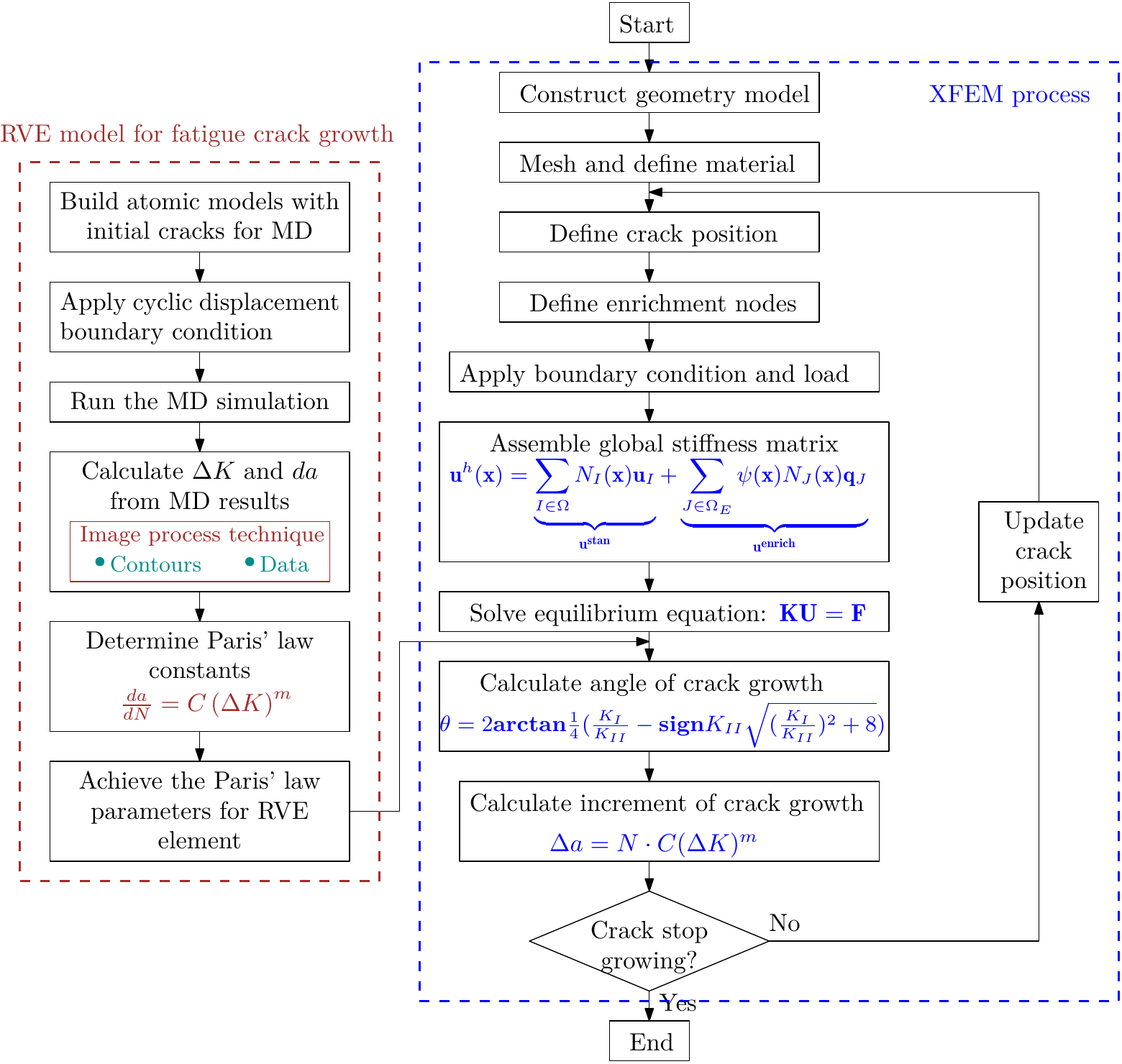}
	\caption{Framework of the RVE based multi-scale method for fatigue crack propagation.}
	\label{fig:flowchart}
	\end{figure*}

	The RVE based multi-scale method is proposed to simulate the fatigue crack propagation on both micro-scale and macro-scale. As shown in Fig. \ref{fig:flowchart}, the RVE based multi-scale method can be divided into two parts: the XFEM process and the RVE model for fatigue crack. The former is used to simulate the fatigue crack propagation on the macro-scale where the XFEM is used here. The latter is used to calculate the Paris constants for fatigue crack propagation by MD simulations on the micro-scale. The Paris law is a bridge between XFEM and MD simulations which connects microscopic property to macroscopic phenomena. Moreover, the image processing technique is used to extract the crack length from MD contours directly. More details are given as follows.

	\subsection{Illustration of the multi-scale model}

\begin{figure*}[htbp] 
		\centering
		\includegraphics[width=\linewidth]{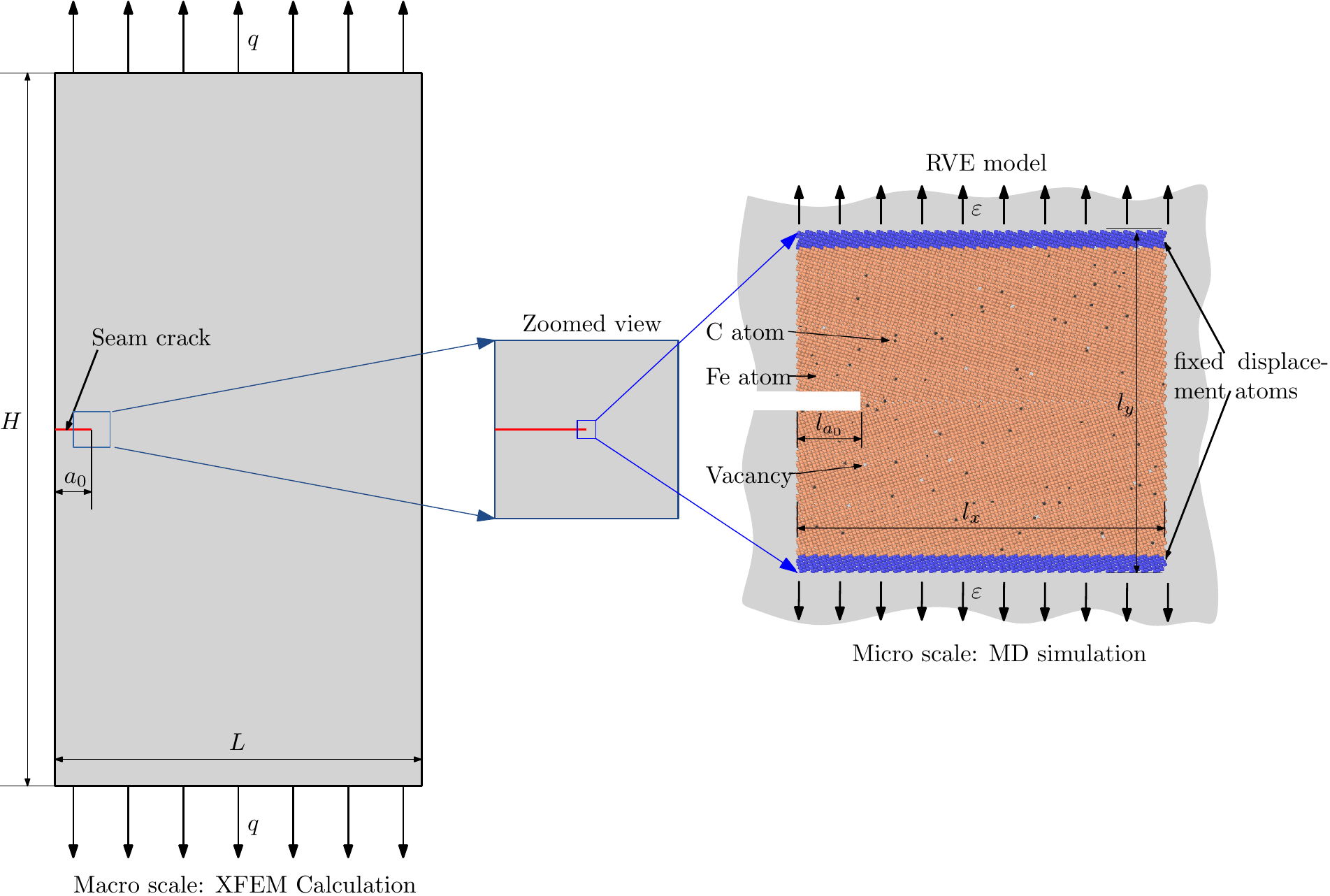}
		\caption{Illustration of the multi-scale model for fatigue crack.}
		\label{fig:model}
		\end{figure*}

	Figure \ref{fig:model} shows a multi-scale model for the fatigue crack propagation under the mode \MakeUppercase{\romannumeral1} loading conditions. The left figure shows the macroscopic model  where a rectangular plate with size of $L \times H$ is considered. A crack with length of $a_0$ is set at the center of the plate.  The specimen is made of conventional mild carbon steel and the material parameters can be found from the study of Bo{\v{z}}i{\'{c}} et al.\cite{Bozic2011}. The geometry and material parameters are listed  in Tab.\ref{tab:material}, where $E, \mu$ and $G$ represent Young modulus, Poisson ratio and shear modulus respectively. The symbol $\sigma_c$ means the yield stress of the specimen.  
	
		\begin{table*}[b]
		\centering
		\small
		\caption{Geometry and material parameters for the macroscopic model}
		\label{tab:material}
		\begin{tabular}{c c c c c c c c}
			\toprule
			$L(mm)$ & $H(mm)$ & $a_0(mm)$ & $q(N/mm) $ & $E(Gpa)$ & $\mu$ & $G(GPa)$ & $\sigma_c(MPa) $ \\ \midrule
			  60    & 120     & 10        & 50        & 206      & 0.3   & 80       & 235              \\ \bottomrule
		\end{tabular}
		\end{table*}
	
	The microscopic model is placed right at the tip of the macroscopic crack as shown in the right of Fig. \ref{fig:model}. The geometry of the atomic model is $200\AA \times 200\AA \times 10\AA$, which has about 34000 atoms in the system. As mentioned above, the specimen is made by conventional mild carbon steel, so the atomic model is made up of iron(Fe) and carbon (C) atoms in which Fe atoms are arranged with body center cubic (BCC) crystal structure and C atoms are inserted as substitutional and interstitial atoms. The ratio of C atoms is $0.2\%$ and the vacancies are also considered in this study. The crystal is in the cubic orientation with $X=[100], Y=[010], Z=[001]$ and the lattice constant of Fe is $2.85 \AA$. As shown in Fig. \ref{fig:model}, a crack is placed on the left edge of the model. The length of the initial crack is equal to $40\AA$.
	
\begin{figure*}[htbp] 
	 	\centering
	 	\includegraphics[]{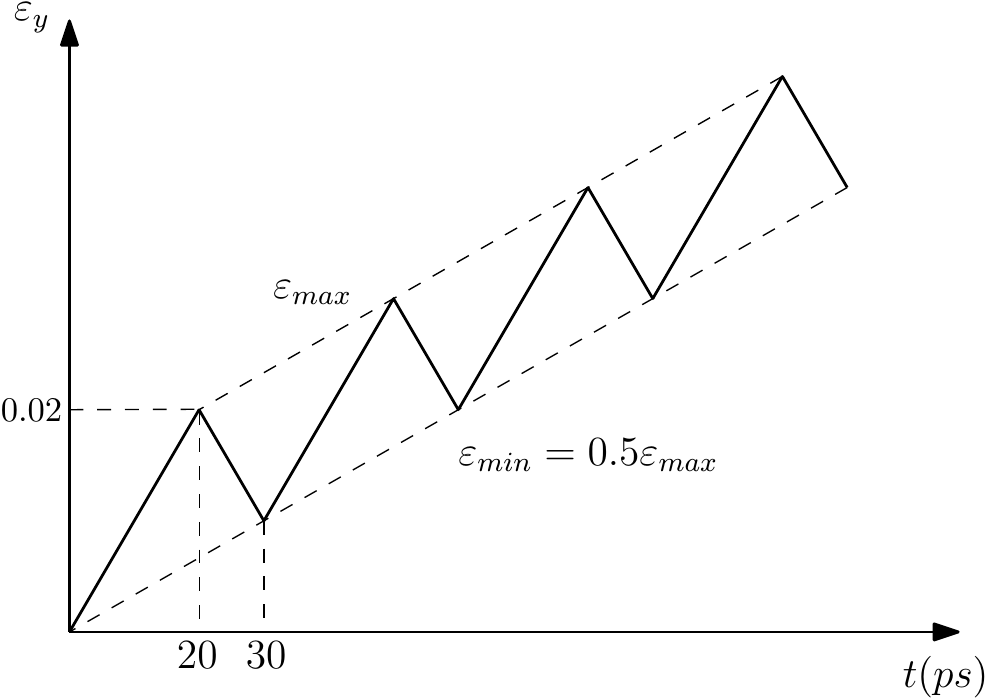}
	 	\caption{Cyclic loading curve with increasing strain amplitude loading.}
	 	\label{fig:rate}
	 \end{figure*}

	 In order to apply the mode \MakeUppercase{\romannumeral1} loading conditions, atoms in the top and bottom layers are fixed as boundaries. Usually, the simulation process of MD includes two steps: relaxation and loading. The relaxation is used to minimize the energy to a equilibrium state. The loading step is to apply the boundary conditions. In this study, a strain-controlled cyclic loading is applied to the fixed atoms labeled in Fig. \ref{fig:model} and  the constant strain rate $r_\varepsilon$  was $1 \times 10^9$. An increasing strain amplitude loading with the load ratio of 0.5 ($R=\varepsilon_{min}/ \varepsilon_{max}$) is used and the cyclic loading curve is shown in Fig. \ref{fig:rate}.
	 
	  For MD simulation, the time-step of 0.001 ps (picoseconds) is used in this study. An open source MD code, named Large-scale Atomic Molecular Massively Parallel Simulator (LAMMPS) \cite{plimpton1995fast}, is used to execute the calculation of MD simulations and the visualization of results is processed by Open Visualization Tool (OVITO) \cite{ovito}.

	\subsection{Determine Paris law constants}

	As mentioned before, the Paris law is the bridge between XFEM and MD simulations and the RVE model built by atoms is used to fit the Eq.\eqref{eq:1} and obtain Paris law constants. According to the Eq.\eqref{eq:1}, it can be found that the stress intensity factor $K_I$ and the crack length $a$ during each loading cycle are essential variables for fitting the $da/dN$ versus $\Delta K$ curve (Fig.\ref{fig:paris}). Therefore, how to get values of  $K_I$ and $a$ is the first step needs to be completed. 
	
	As for the stress intensity factor $K_I$, since the LAMMPS can calculate the atomic stress by Eq.\eqref{eq:6}, so the $K_I$ can be easily calculated by the Griffith level by Eq.\eqref{eq:8}, where $\sigma_y$ means the stress component of crack tip \cite{ma2014molecular}, and $a$ is the crack length.  Then the range of stress intensity factor $\Delta K$ can be calculated by Eq.\eqref{eq:9}, where values of the maximum ($\sigma_{max}$) and minimum ($\sigma_{min}$) stresses are used to calculate $K_{max}$ and $K_{min}$ by Eq.\eqref{eq:8}.
	
	\begin{gather}
		K_I=\sigma_y \sqrt{\pi a}, \label{eq:8}\\
		\Delta K=K_{max}-K_{min}. \label{eq:9}
	\end{gather}
	
\begin{figure*}[htbp] 
		\centering
		\includegraphics[width=\linewidth]{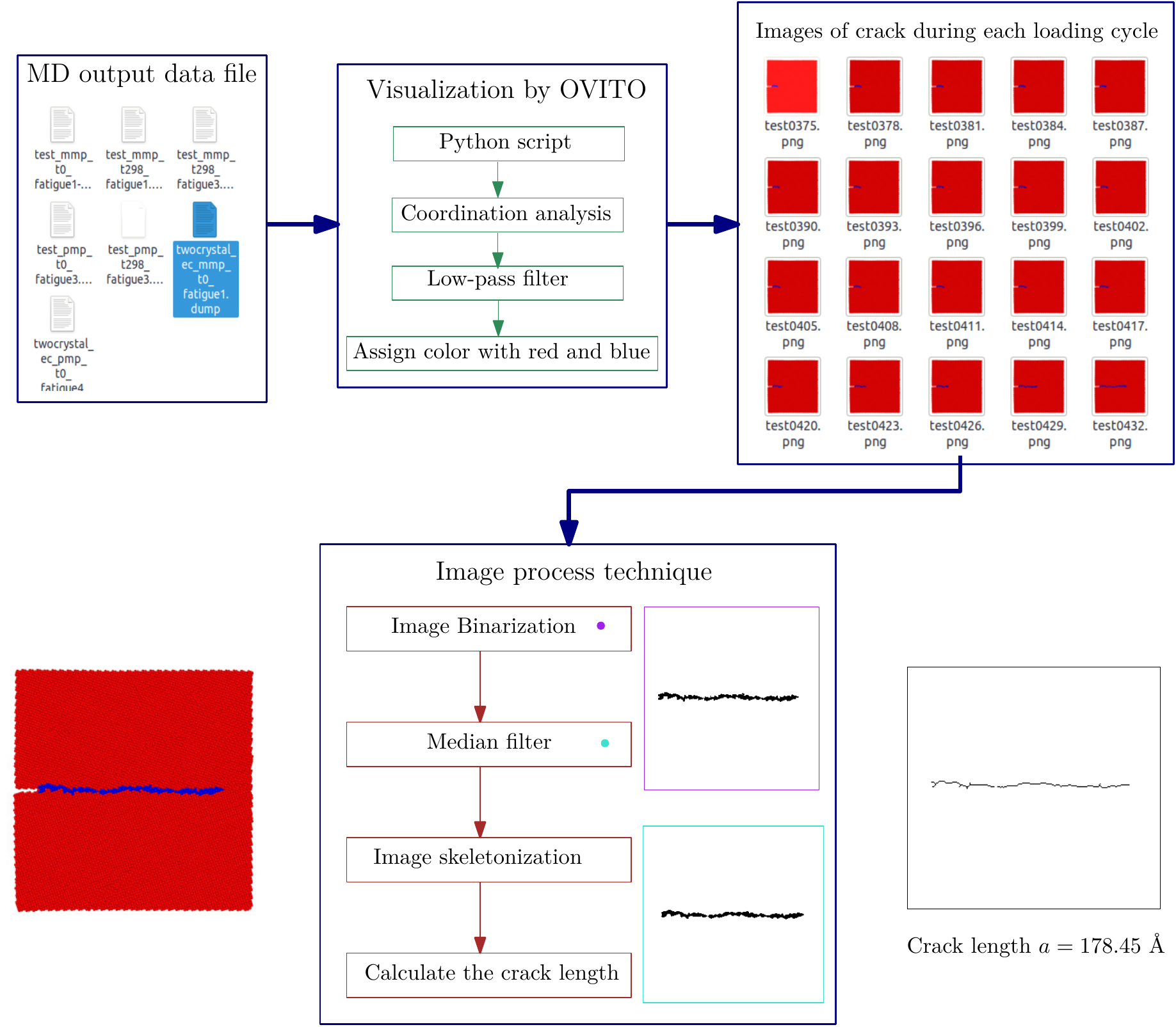}
		\caption{Illustration of the image based crack extracting method.}
		\label{fig:imgpt}
	\end{figure*}

	As for the crack length $a$, it is difficult to calculate it from  MD result files because the LAMMPS cannot recognize the position of crack tips in real time and the size of  MD results files is usually too huge to read. Therefore, an image based crack extracting method is suggested to extract the crack path by image processing technique in this study. Considering the huge memory requirement of MD result files, the main superiority of the image based method is to reduce the memory requirement and improve the efficiency of crack extracting process by using a contour image (several hundred kilobytes) as the input instead of an MD result file (several hundred gigabytes). The illustration of the image based crack extracting method is shown in Fig. \ref{fig:imgpt}. As mentioned above, the visualization of  MD result files is processed by OVITO, in which a python script is written to generate corresponding contours automatically. The coordination analysis counts the number of atoms for each atom that are within the cutoff range around its position, so the boundaries and fractures can be recognized with a low-pass filter because less neighbor atoms are around its position than others. Then atoms around crack will be assigned to blue color after removing the boundary atoms. After that, the image process technique is used to extract the skeleton of crack with operations of binarization, median filtering and skeletonization. Finally, the crack length can be calculated out from the skeleton of crack.
	
	After getting the data of the range of stress intensity factor ($\Delta K$) and the increment of crack ($da/dN$), the Paris law constants $C$ and $m$ can be calculated by fitting the Eq.\eqref{eq:1} with the data. It is necessary to introduce experimental Paris law constants for the material of specimen to provide a validation data for the next numerical results, which gives $m=2.75$ and $C=1.43 \times 10^{-11}$ \cite{Bozic2011}.

	\subsection{XFEM for fatigue crack propagation}
	Extended finite element method,  a popular numerical method for simulating crack propagation \cite{Moes1999,Belytschko2009}, is used to investigate the fatigue crack propagation in this study. The XFEM was an improved method based on the conventional finite element method by adding enrichment functions into the shape function to describe the discontinuities. Generally, the displacement approximation of XFEM can be described as Eq.\eqref{eq:10}, where $N_I$ and $\bm u_I$ represent the standard FEM shape function and nodal degrees of freedom (DOF) respectively. The $\psi (\bm x)$ is the enrichment function and the $\bm q_J$ is the additional nodal degree of freedom. Rewrite the Eq. \eqref{eq:10} as Eq.\eqref{eq:11}, where $\Omega$ is the solution domain, $\Omega_s$ is the domain cut by crack,  $\Omega_{T}$ is the domain which crack tip located, $H({\bm x})$ is the shifted Heaviside enrichment and  $\Phi_{\alpha}(\bm{x})$ is the shifted crack tip enrichment. The details of $H({\bm x})$ and $\Phi_{\alpha}(\bm{x})$ are defined as Eqs. \eqref{eq:12} and \eqref{eq:13}. Then the discrete equilibrium equation can be obtained by the principle of virtual work. Finally, the displacement result can be obtained after solving the equilibrium equation.
	
	\begin{gather}
		\bm u^h (\bm x )=\underbrace{\sum_{I \in \Omega} N_I (\bm x) \mathbf u_I}_{\bm u^\textbf{stan}} + \underbrace{\sum_{J \in \Omega_E} \psi (\bm x) N_J (\bm x) \bm q_J}_{\bm u^\textbf{enrich}}, \label{eq:10} \\ 
		\bm{u}^{h}(\bm{x})=\sum_{I \in \Omega} N_{I}(\bm{x}) \bm{u}_{I}+\sum_{I \in \Omega_{S}} H_{I}(\bm{x}) N_{I}(\bm{x}) \bm{a}_{I}+\sum_{I \in \Omega_{T}} \sum_{\alpha=1}^{4} \Phi_{I, \alpha}(\bm{x}) N_{I}(\bm{x}) \bm{b}_{I}^{\alpha}, \label{eq:11} \\
		H(\bm{x})=\left\{
		\begin{array}{ll}
				+1 & \text { Above crack } \\
				-1 & \text { Below crack }
		\end{array}\right., \label{eq:12}\\
		\left\{\Phi_{\alpha}(\bm{x})\right\}_{\alpha=1}^{4}=\sqrt{\bm{r}}\left\{\sin \frac{\theta}{2}, \cos \frac{\theta}{2}, \sin \theta \sin \frac{\theta}{2}, \sin \theta \cos \frac{\theta}{2}\right\}. \label{eq:13}
	\end{gather}

	As for the definition of crack propagation, the direction and magnitude of crack propagation are usually used to determine how the crack will grow. The direction can be calculated by the maximum circumferential stress criterion \cite{erdogan1963crack}. The angle of crack propagation is calculated by
	
	\begin{equation}
		\theta=2 \arctan \frac{1}{4} (\frac{K_I}{K_{II}}-\textbf{sign} K_{II} \sqrt{(\frac{K_I}{K_{II}})^2+8}),
	\end{equation}

	where $\theta$ is defined in the crack tip coordinate system, $K_I$ and $K_{II}$ are the mixed-mode stress intensity factors. The details are given in the reference \cite{erdogan1963crack}. For the magnitude of crack propagation, especially for the fatigue crack propagation, the increment of crack propagation $\Delta a$ can be calculated by the Paris law which can be write as Eq.\eqref{eq:15}, where the constants $m$ and $C$ are obtained from the RVE model described above.

	\begin{equation} \label{eq:15}
				\Delta a = N \cdot C (\Delta K)^m
	\end{equation}

	\section{Results and discussions}
	
	\subsection{Paris law constants calculation}
	
	\begin{figure*}[htbp]
		\centering
		\subfigure[]{\includegraphics[width=2in]{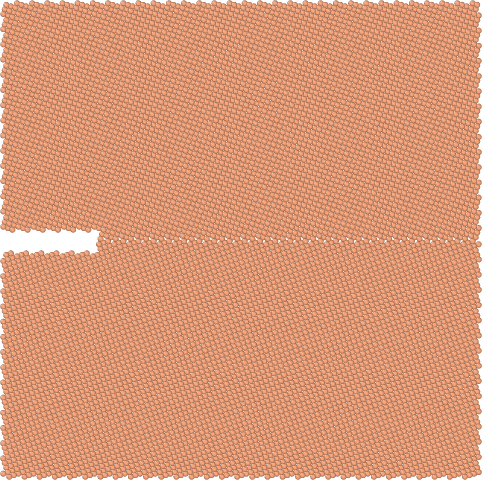}}
		\hspace{0.1in}
		\subfigure[]{\includegraphics[width=2in]{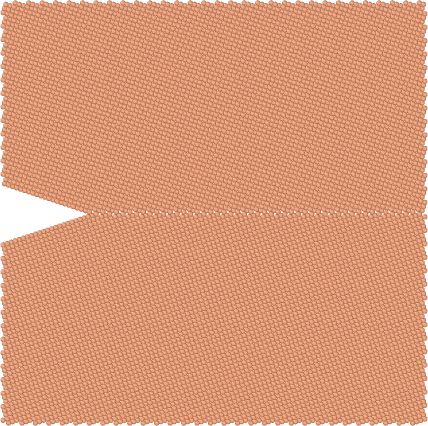}}\\
		
		\subfigure[]{\includegraphics[width=2in]{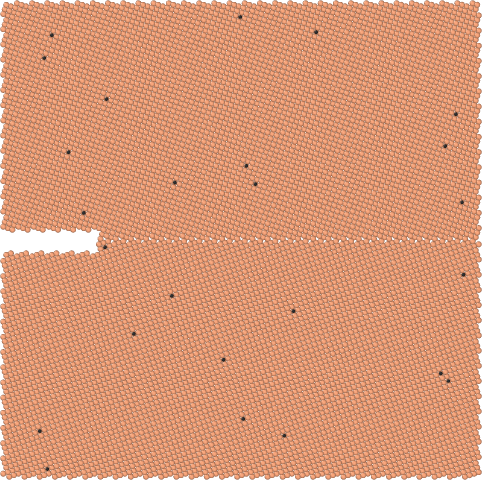}}
		\hspace{0.1in}
		\subfigure[]{\includegraphics[width=2in]{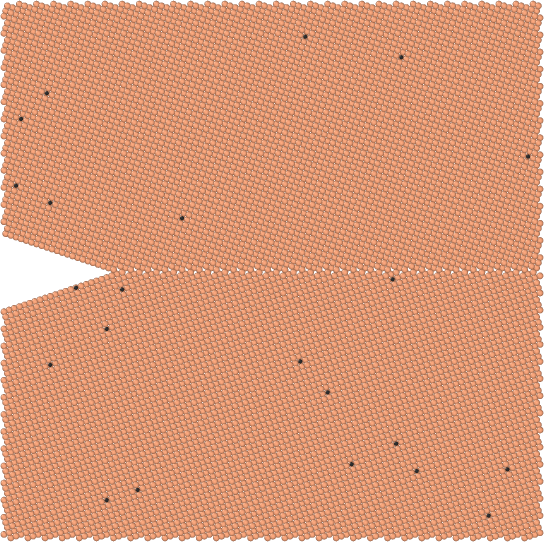}}
		\caption{Atomic model with, (a) \& (c)blunt edge crack, (b) \& (d) sharp edge crack,  which made by, (a) \& (b) Pure iron atoms, (c) \& (d) Iron atoms with 0.2\% carbon atoms and 0.5\% vacancies. }
		\label{fig:atomicmodels}
	\end{figure*}
	
\begin{figure*}[htbp]
		\centering
		\subfigure[]{\includegraphics[width=3in]{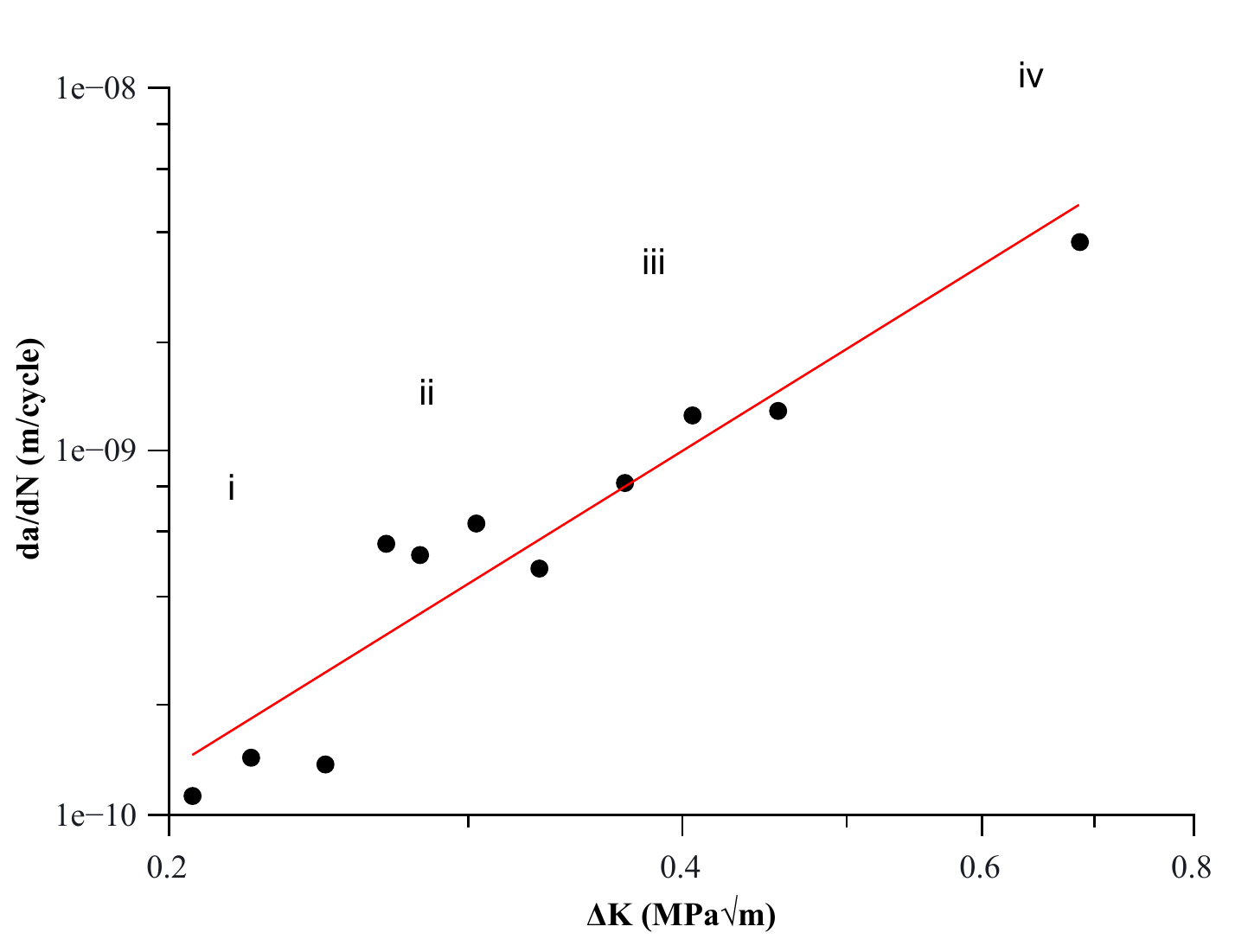}}
		\hspace{0.1in}
		\subfigure[]{\includegraphics[width=2in]{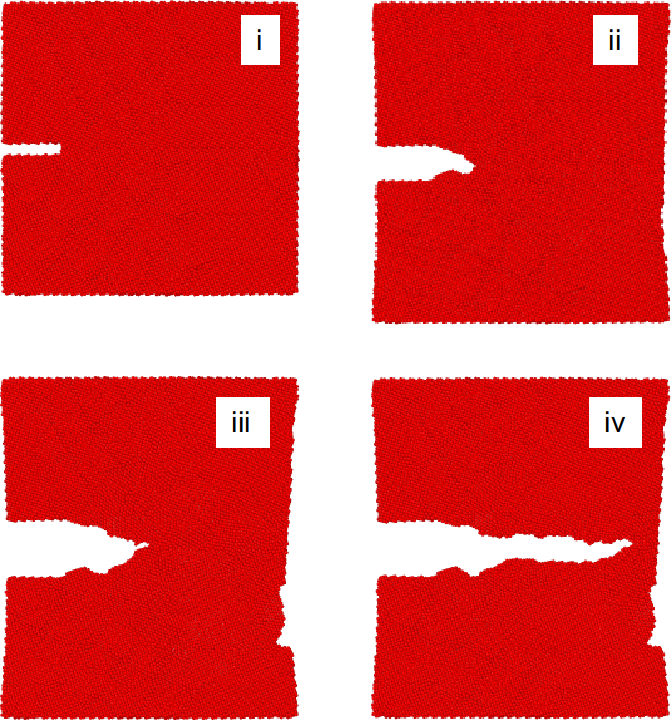}}\\
		
		\subfigure[]{\includegraphics[width=3in]{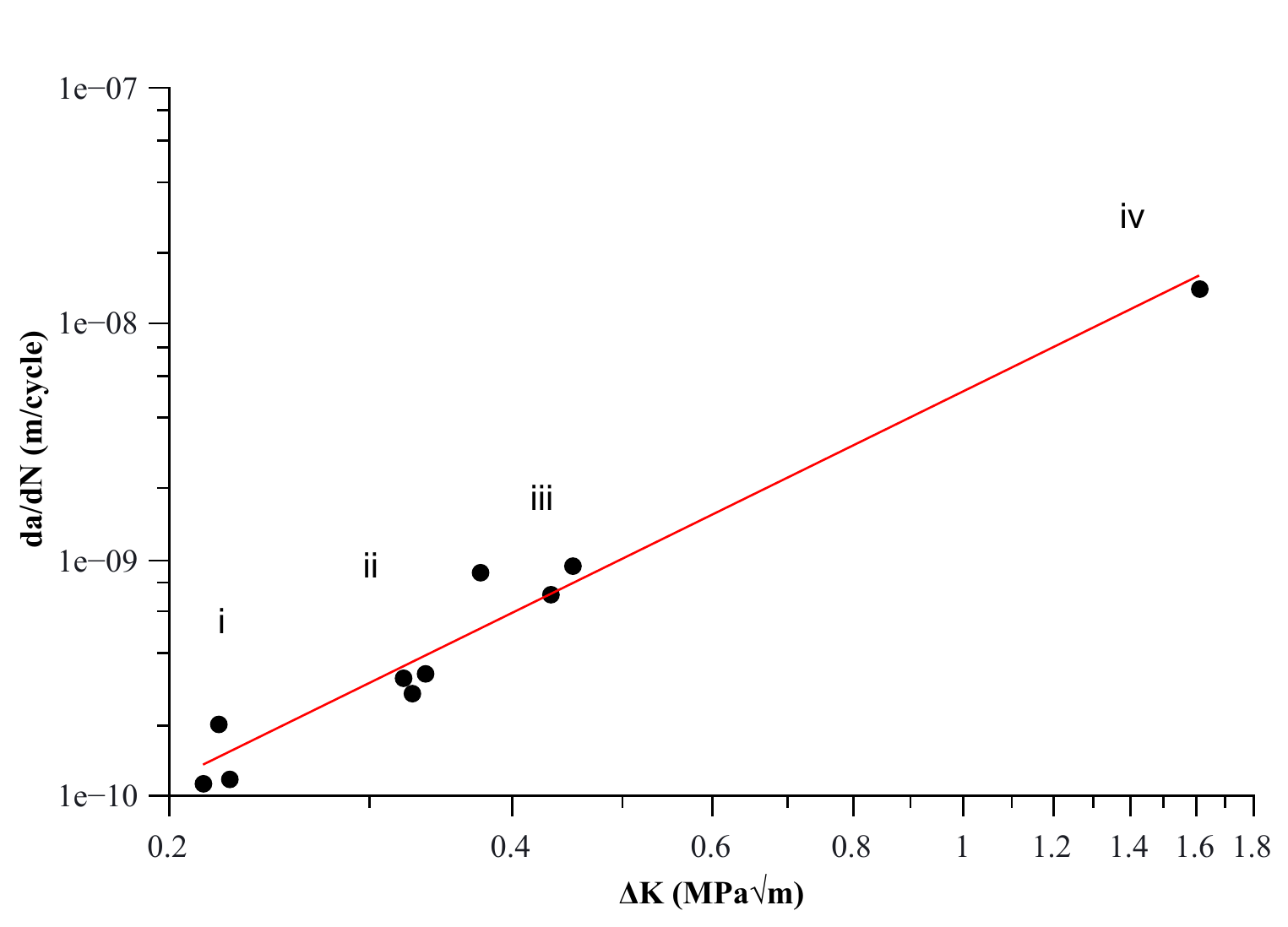}}
		\hspace{0.1in}
		\subfigure[]{\includegraphics[width=2in]{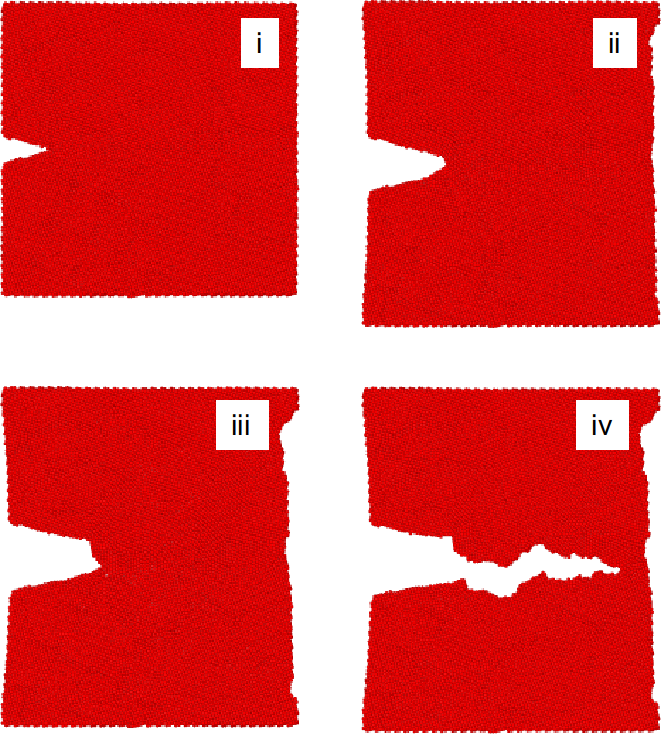}}
		
		\caption{Fatigue growth rate curve derived by MD models with, (a) blunt edge crack, and (c) sharp edge crack,  which made by pure iron atoms.}
		\label{fig:parisresults}
	\end{figure*}

\begin{figure*}[htpp]
		\centering
		\subfigure[]{\includegraphics[width=3in]{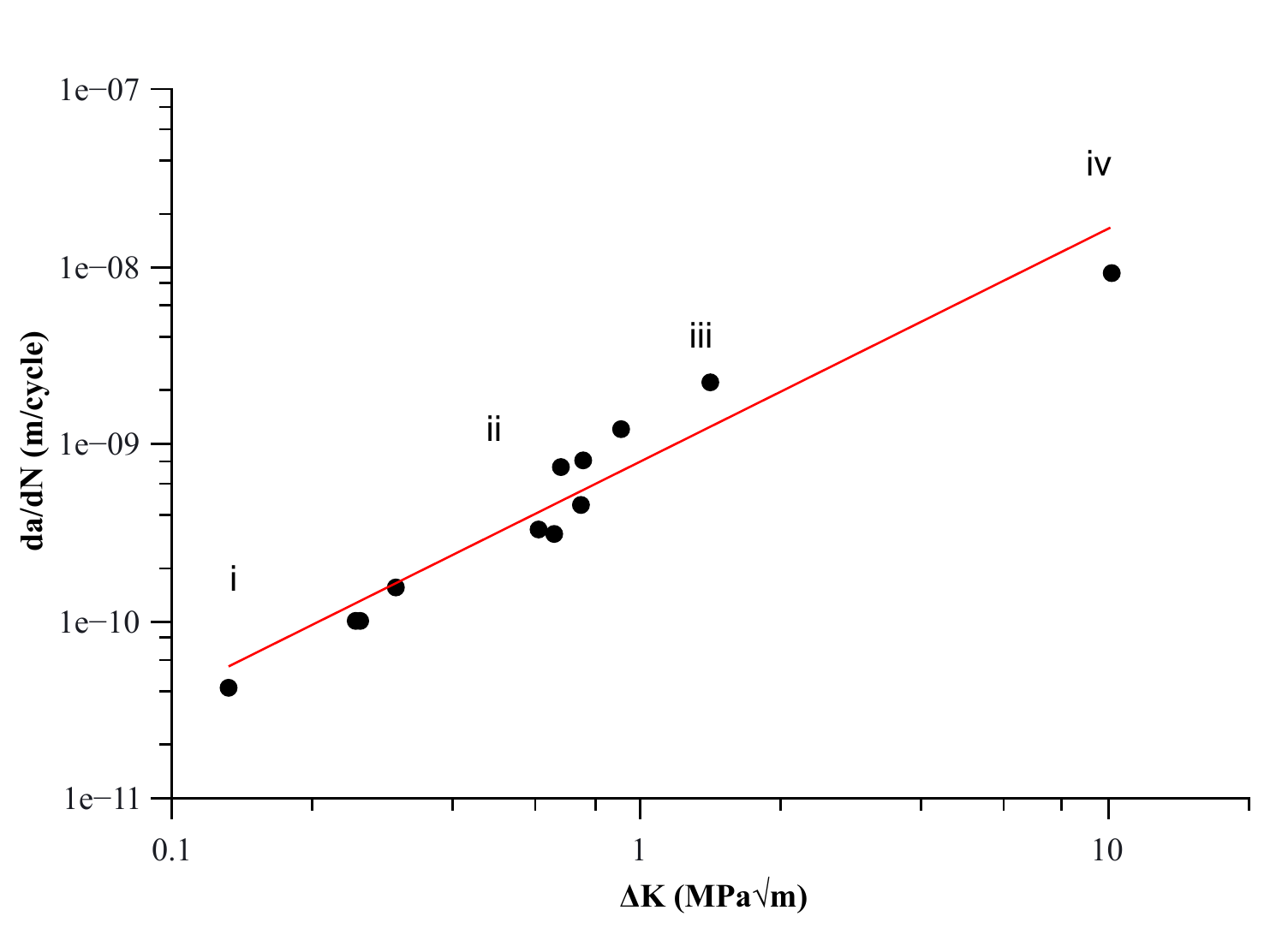}}
		\hspace{0.1in}
		\subfigure[]{\includegraphics[width=2in]{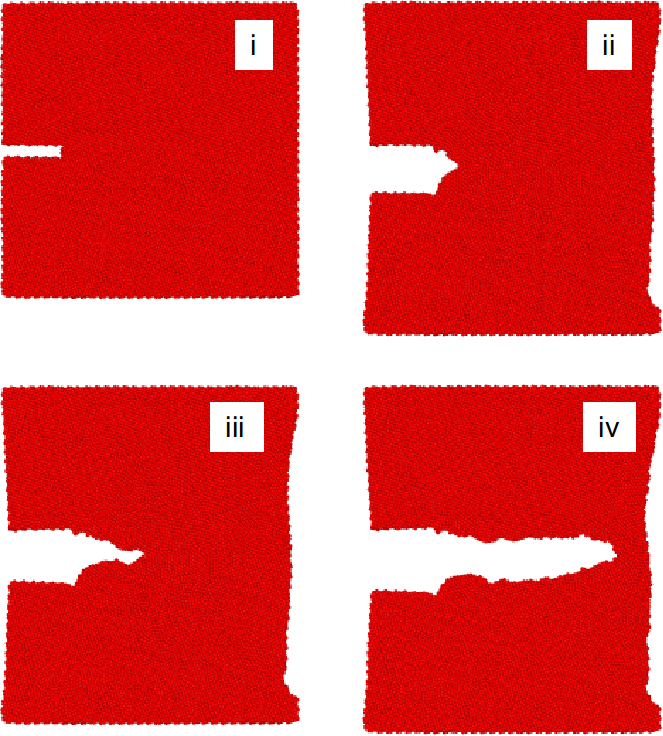}}\\
		
		\subfigure[]{\includegraphics[width=3in]{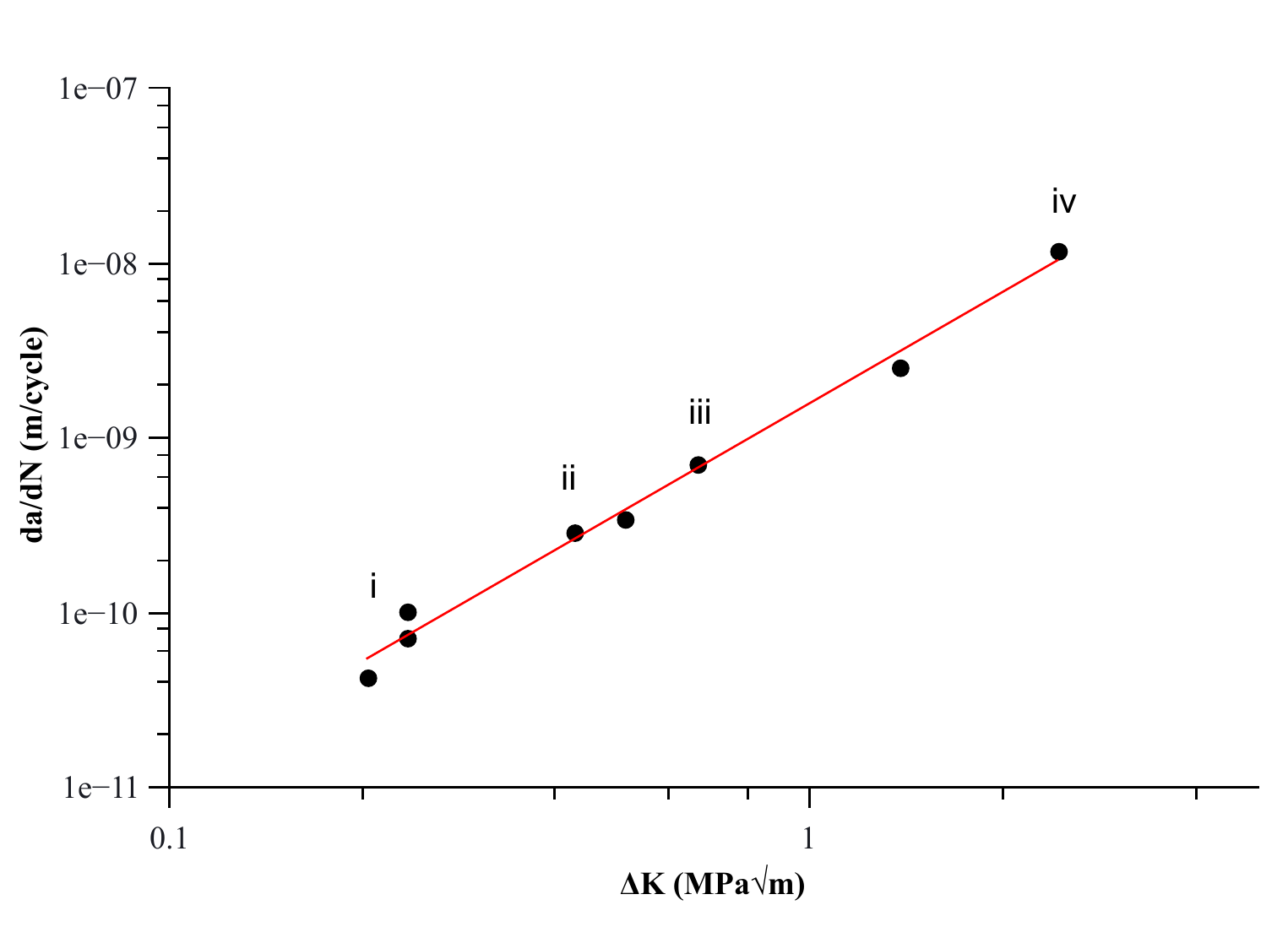}}
		\hspace{0.1in}
		\subfigure[]{\includegraphics[width=2in]{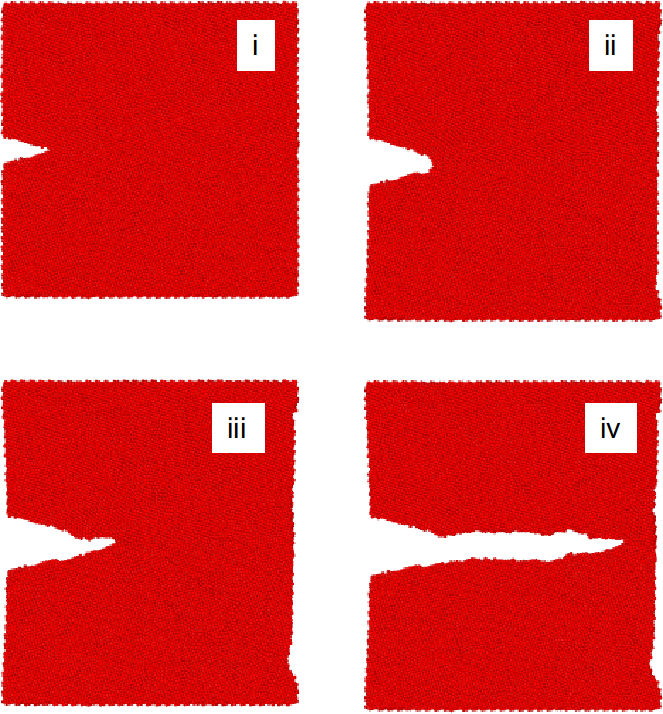}}
		
		\caption{Fatigue growth rate curve derived by MD models with, (a) blunt edge crack, and (c) sharp edge crack, which made by iron atoms with 0.2\% carbon atoms and 0.5\% vacancies. }
		\label{fig:parisresults2}
	\end{figure*}

	To fully investigate the behavior of fatigue crack propagation, four atomic models with various micro-defect patterns for each type are constructed to execute by LAMMPS, which are characterized  by  Fig. \ref{fig:atomicmodels}, where (a) is the model with blunt edge crack which made by pure iron atoms, (b) is the model with sharp edge crack which made by pure iron atoms. Figure \ref{fig:atomicmodels}(c) and (d) are trying to model the micro-structure of mild carbon steel by iron atoms with 0.2\% carbon atoms and 0.5\% vacancies. Then the strain-controlled cyclic loading with the strain rate of $1 \times 10^9$ (as shown in Fig. \ref{fig:rate}) should be applied to those four models to simulate the process fatigue crack propagation. The fatigue growth rate curves are shown in Fig. \ref{fig:parisresults} and \ref{fig:parisresults2}. Figure \ref{fig:parisresults}(a) and (b) exhibit a progressive cracking behavior for the pure iron model with a blunt edge crack (Model A). Point \romannumeral1\ indicates the initial state. After that the crack starts to grow gradually (point \romannumeral2). As the fatigue crack keeps growing, then the crack comes to state \romannumeral3\ and finally reaches the fracture state (point \romannumeral4). As observed in Fig. \ref{fig:parisresults}(c) and (d), a pure iron plate with a sharp edge crack (Model B) was simulated by MD method. Points \romannumeral1 - \romannumeral4\ indicate the same states as in Fig. \ref{fig:parisresults}  (a)  to make a comparison with model A. Figure \ref{fig:parisresults2} demonstrated the fatigue cracking behaviors in mild carbon steel by considering the interstitial atoms (C) and vacancies, where the plate was made by iron atoms with $0.2\%$ carbon atoms and $0.5\%$ vacancies. Figure \ref{fig:parisresults2}(a) and (b) exhibit the cracking behavior for the mild carbon steel plate with a blunt edge crack while the Fig. \ref{fig:parisresults2} (c) and (d) exhibit the cracking behavior with a sharp edge crack.
	
	As shown in the figures, the Table \ref{tab:constants} shows the values of Paris law constants which are evaluated by fitting the fatigue growth rate curve (Fig. \ref{fig:parisresults} and \ref{fig:parisresults2}) with Paris model. Compared with the experimental Paris law constants in the reference \cite{Bozic2011} ($m=2.75$ and $C=1.43 \times 10^{-11}$ ), it can be found that the result of pure iron model with the blunt edge crack agrees quite well with the experimentally determined values of Paris law constants. Moreover, it is evident that different models and initial crack types result in different calculated value of Paris law constant as shown in Tab. \ref{tab:constants}, i.e., the fatigue crack properties should largely depend on the micro-structure of materials.

\begin{table*}[htbp]
	\centering
	\small
	\caption{Paris law constants derived from MD simulations.}
	\label{tab:constants}
	\begin{tabular}{c c c c c}
		\toprule
		Model number & Material type                                          & Crack type  & $m$      & $C$                      \\ \midrule
		     A       & Pure iron atoms                                        & Blunt crack & $2.9041$ & $1.4299 \times 10^{-11}$ \\
		     B       &                                                        & Sharp crack & $2.3645$ & $5.1902\times 10^{-9}$   \\
		     C       & Iron atoms with 0.2\% carbon atoms and 0.5\% vacancies & Blunt crack & $1.3144$ & $7.9016 \times 10^{-10}$ \\
		     D       &                                                        & Sharp crack & $2.1141$ & $1.5778 \times 10^{-9}$  \\ \bottomrule
	\end{tabular}
\end{table*}

	\subsection{Fatigue crack propagation simulation}
	After getting the Paris law constants by the RVE model, the fatigue crack propagation can be evaluated by the XFEM using those obtained values of $m$ and $C$ on the macro scale. As described in Section 3, the angle of crack growth will be calculated by the maximum circumferential stress criterion and the increment of crack propagation will be determined by the Paris law. Taking model A as an example,	Paris constants obtained from the result of pure iron model with initial blunt crack are used in this section, where the slope  $m=2.9041$ and the y-axis intercept $C=1.4299 \times 10^{-11}$. 
	
\begin{figure*}[htbp] 
	\centering
	\includegraphics[width=5in]{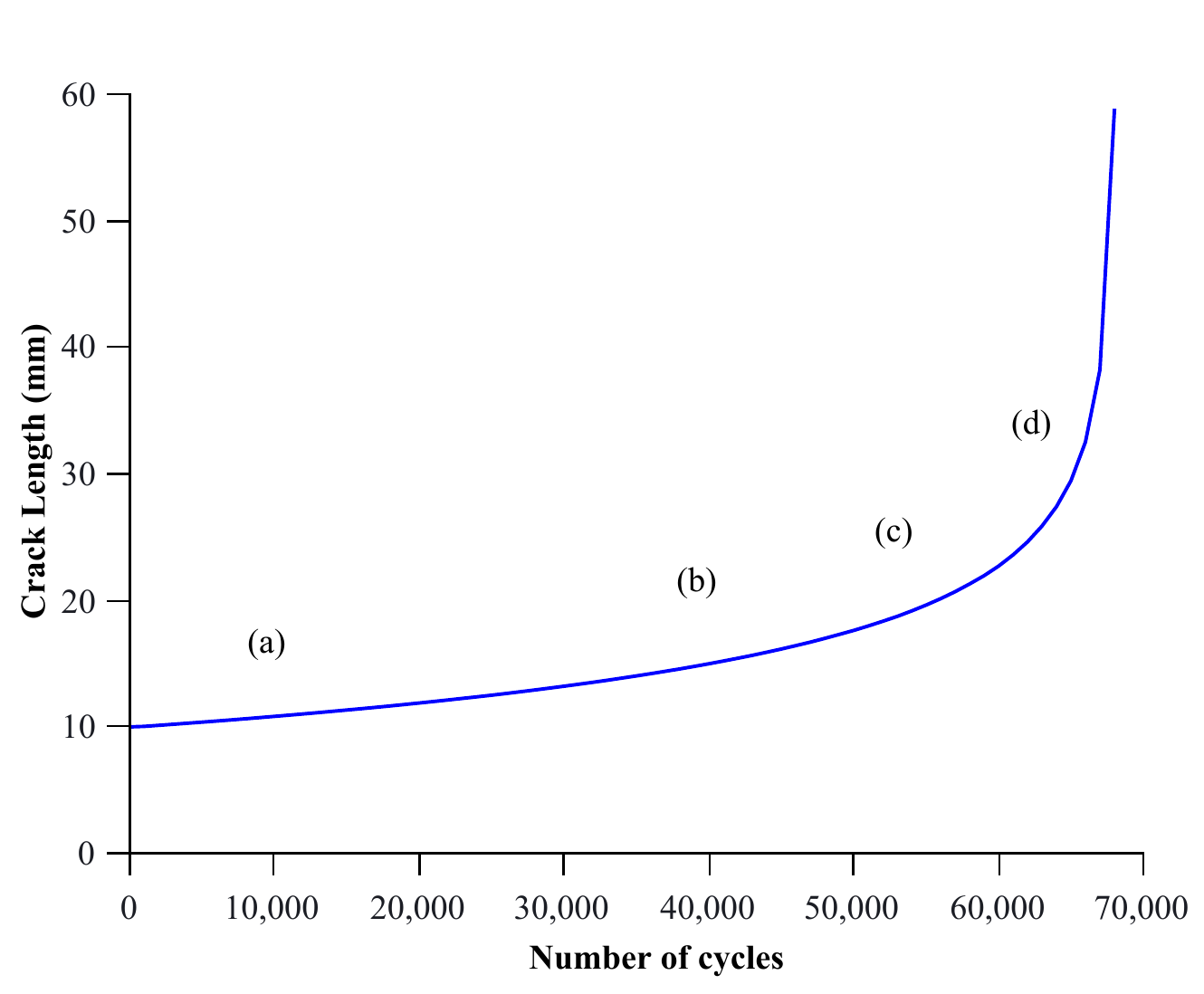}
	\caption{Fatigue crack length versus number of loading cycles.}
	\label{fig:cracklength}
\end{figure*}

\begin{figure*}[htbp]
	\centering
	\subfigure[]{\includegraphics[width=1.5in]{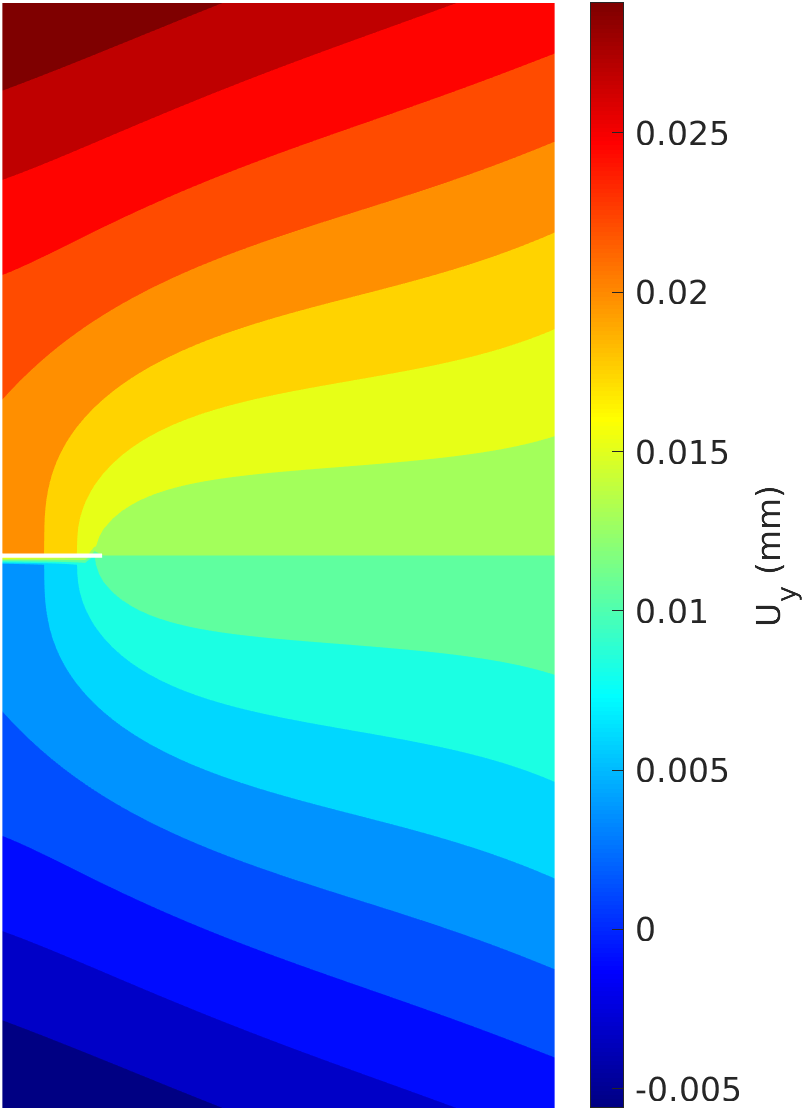}}
	\hspace{0.1in}
	\subfigure[]{\includegraphics[width=1.5in]{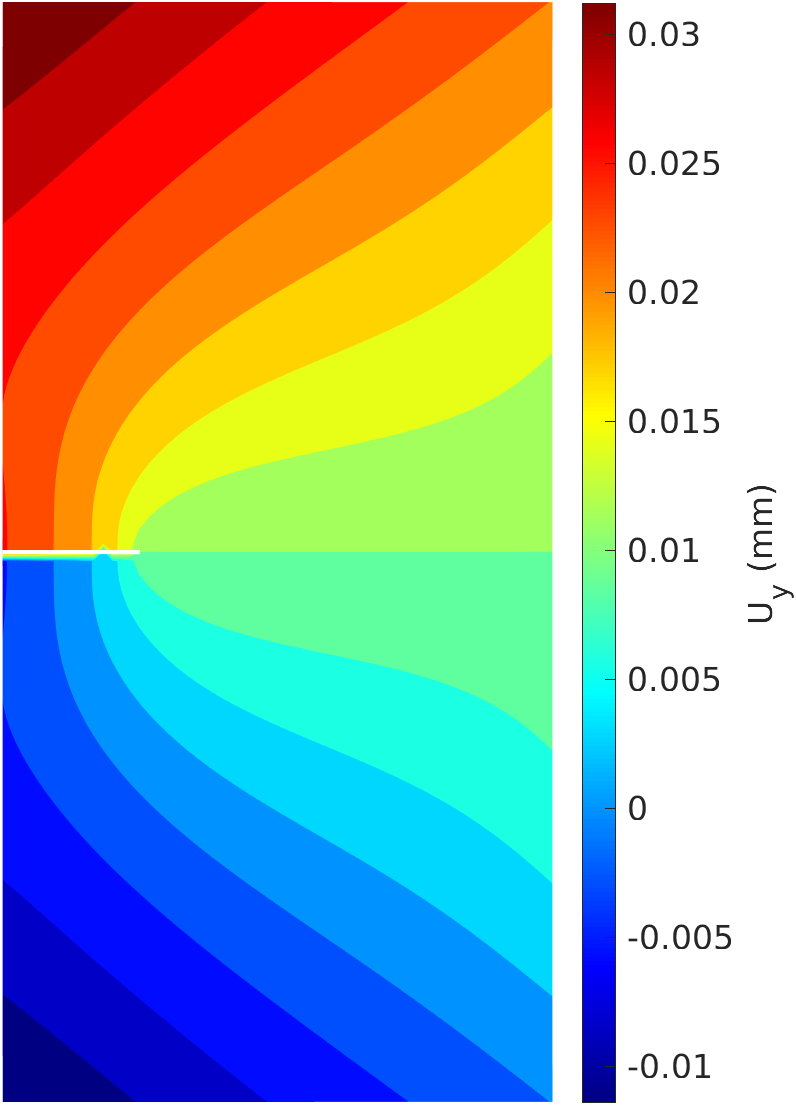}}
	\hspace{0.1in}
	\subfigure[]{\includegraphics[width=1.5in]{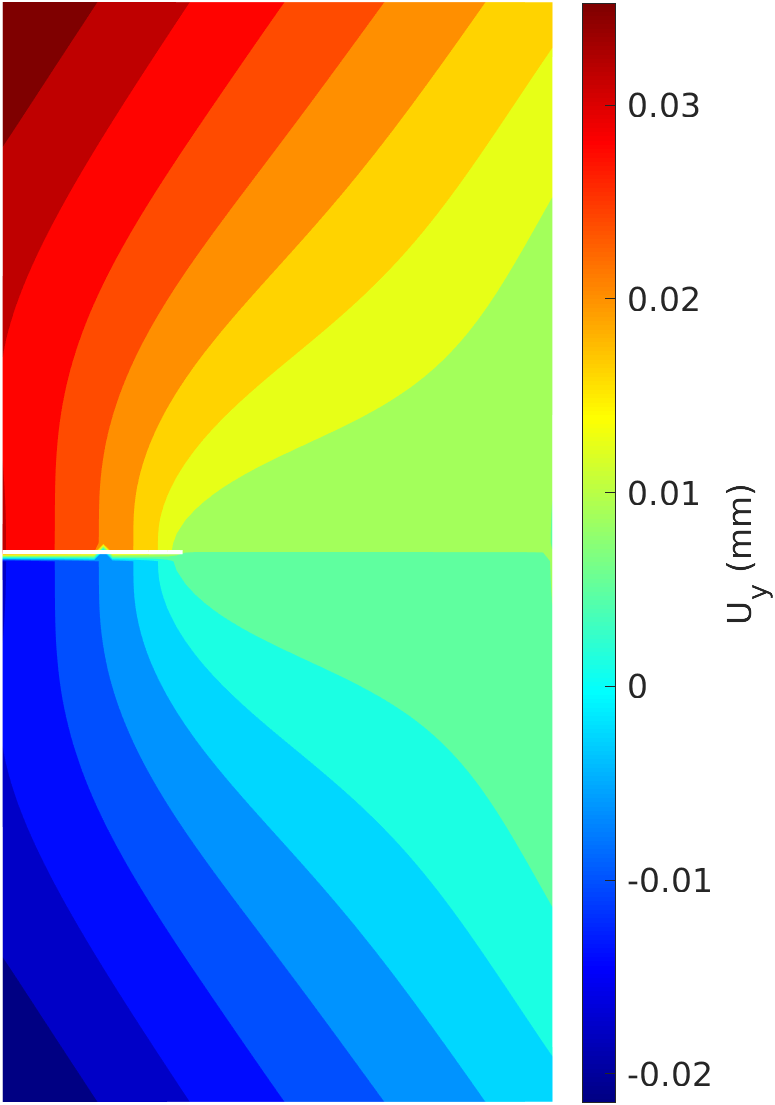}}
	\hspace{0.1in}
	\subfigure[]{\includegraphics[width=1.5in]{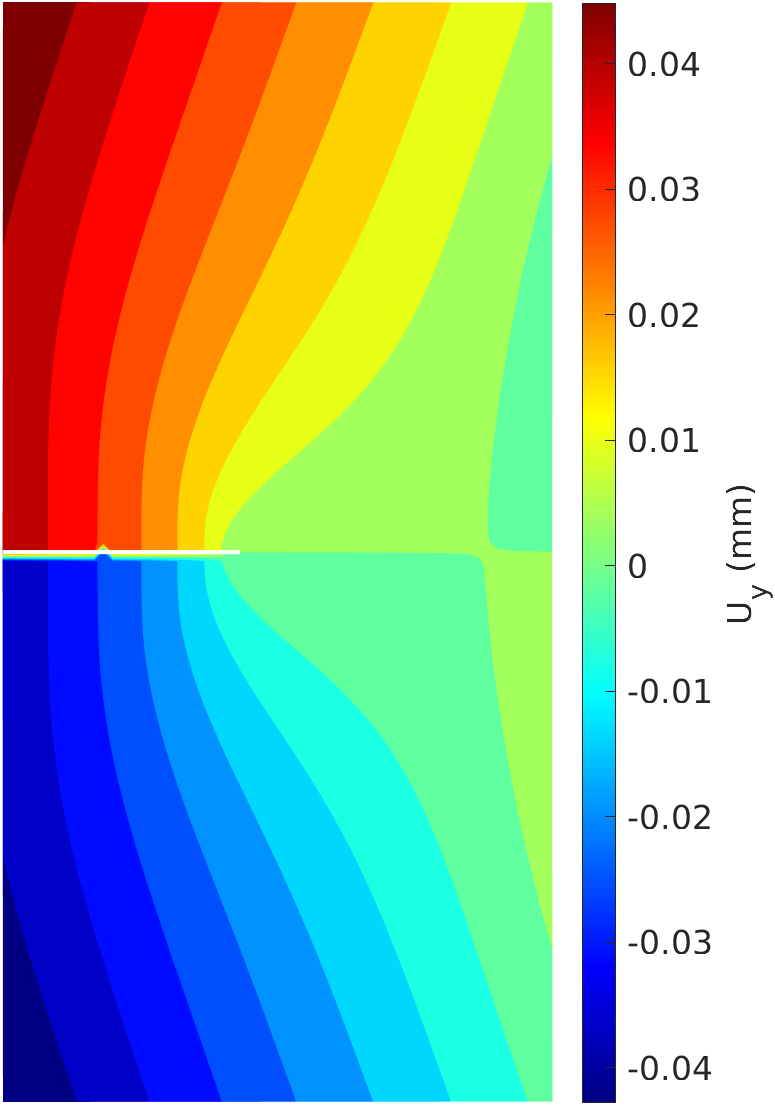}}
	
	\caption{Contours of displacement in y-axis at (a) 10,000th, (b) 40,000th, (c) 55,000th, and (d) 63,000th loading cycles.}
	\label{fig:disp}
\end{figure*}

\begin{figure*}[htbp]
	\centering
	\subfigure[]{\includegraphics[width=1.5in]{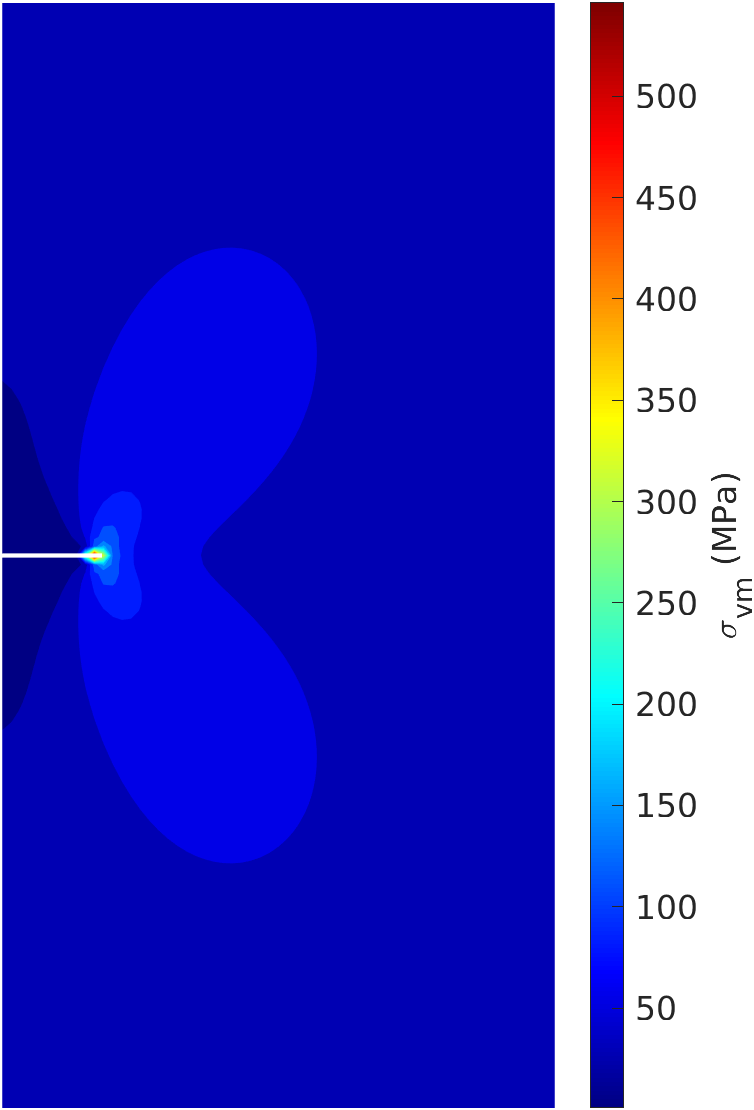}}
	\hspace{0.1in}
	\subfigure[]{\includegraphics[width=1.5in]{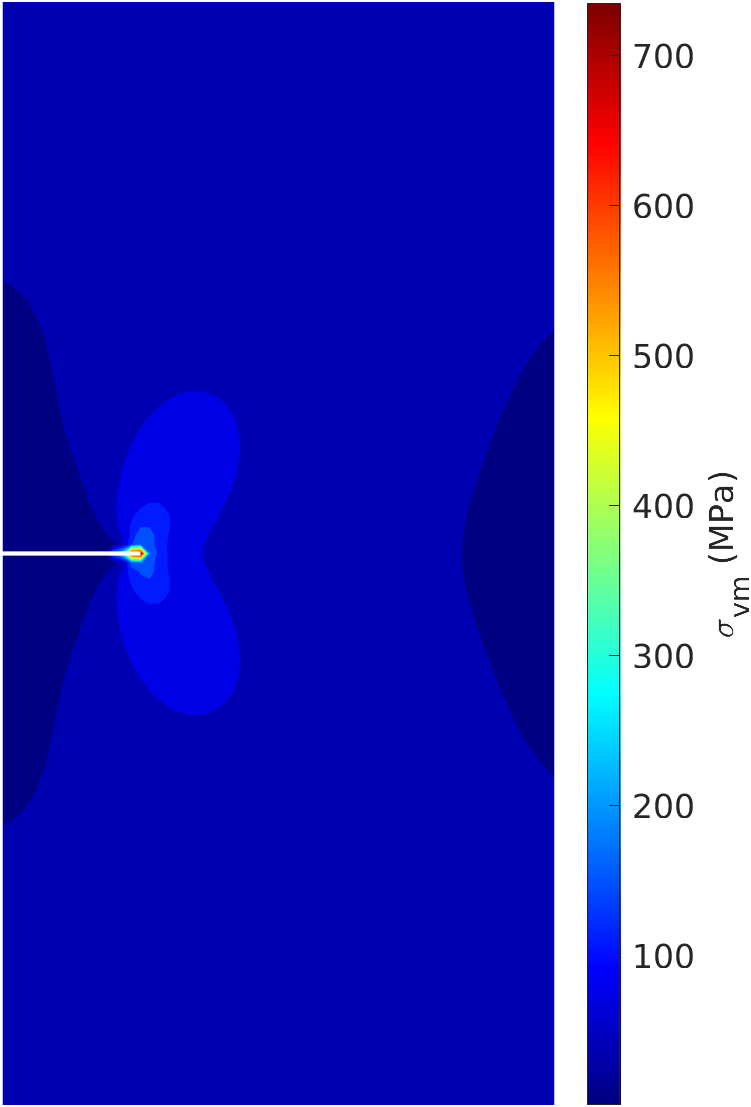}}
	\hspace{0.1in}
	\subfigure[]{\includegraphics[width=1.5in]{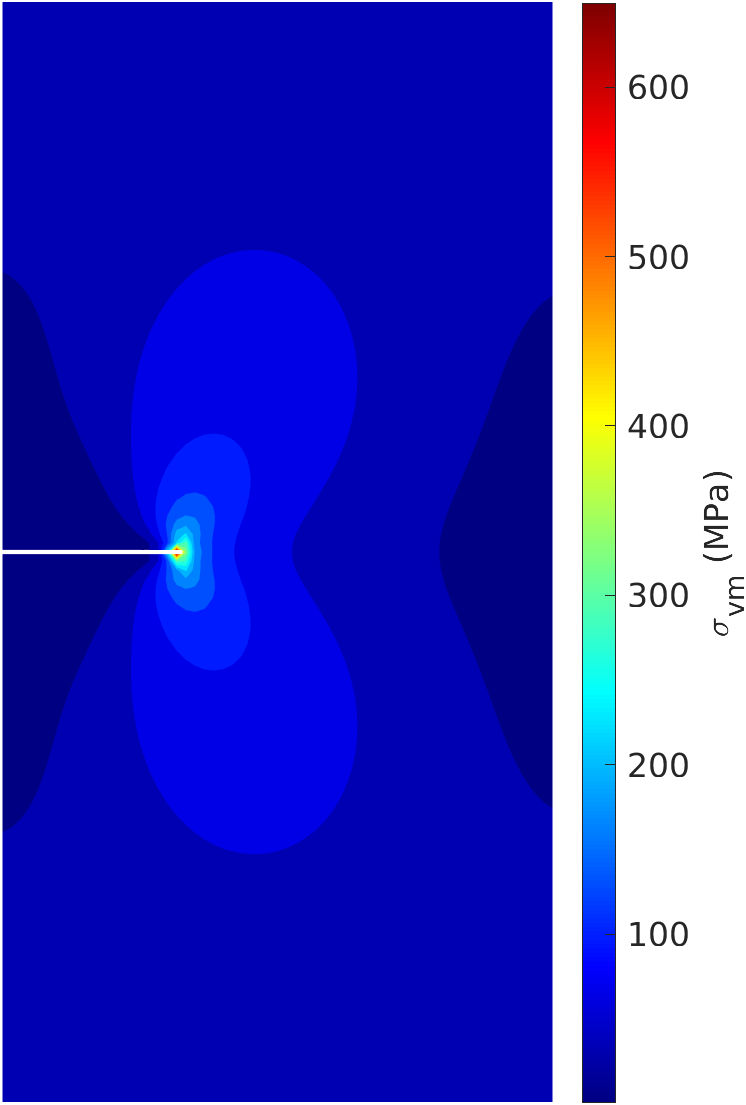}}
	\hspace{0.1in}
	\subfigure[]{\includegraphics[width=1.5in]{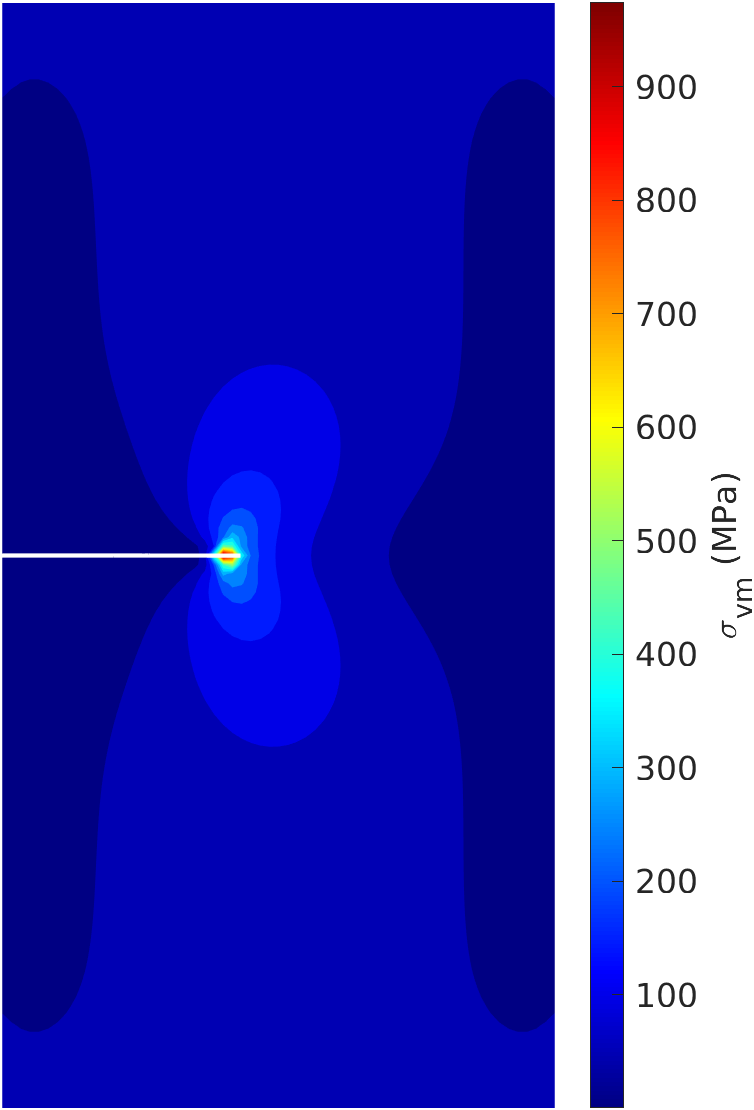}}
	
	\caption{Contours of Von-Mises stress at (a) 10000th, (b) 40000th, (c) 55000th, and (d) 63000th loading cycles.}
	\label{fig:stress}
\end{figure*}

	As shown in Fig. \ref{fig:model}, a rectangular plate with an edge crack is considered in this study. The geometry and material parameters are shown in Tab. \ref{tab:material}. Then the X-FEM is used to simulate the fatigue crack propagation. The curve of fatigue crack length versus number of  loading cycles   is shown in Fig. \ref{fig:cracklength}. Figures \ref{fig:disp} and \ref{fig:stress} demonstrate the contours of displacement and stress at 10,000th, 40,000th, 55,000th and 63,000th loading cycles respectively. It is clearly observed that the fatigue crack length is exponentially growing as the increasing of the number of  loading cycles. According to the curve in Fig. \ref{fig:cracklength}, it can be determined that the plate  will fracture at approximately 65,000th  loading cycles.
	
	\section{conclusions}
	
	In this study,a representative volume element (RVE) based multi-scale method is proposed to study the mechanism of crack propagation in this study, where the extended finite element (XFEM) and MD
	methods are used to study the behavior of macro and micro crack propagation respectively. In order to connect the microscopic property to the macroscopic phenomenon, the Paris law is used as a bridge between XFEM and MD methods. The Paris law constants will be calculated from the result of MD simulation and then the constants will be transferred to XFEM to calculate the life of specimen. The main conclusions were summarized as follows:
	
	\begin{enumerate}[(1)]
		\item The effect of micro-structural defects including interstitial atoms, vacancies have been considered. The results indicate that the micro-structural defects can deeply influence the values of Paris law constants.
		\item Considering the huge memory requirement of MD result files, the image based crack extracting method can significantly reduce the memory requirement and improve the efficiency of crack extracting process.
		\item  The behaviors of fatigue crack  both on micro and macro scale  are investigated by the RVE based multi-scale method and the  life of the specimen can be also evaluated.
	\end{enumerate}
	
	\section*{acknowledgments}
	This work was partially supported by Project of the Key Program of National Natural Science Foundation of China (Grant Number 11972155).
	
	\bibliographystyle{unsrt}

\end{document}